\documentclass[12pt]{article}
\newcommand{\ncm}{\newcommand}
\ncm{\aut}{auto\-mor\-phi\-sm} \ncm{\Inn}{\mbox{\rm Inn($A$)}}
\ncm{\Ap}{\mbox{$\overline{\rm Inn}(A)$}} \ncm{\Ext}{\mbox{\rm
Ext}} \ncm{\Ex}{\mbox{\rm Ex}} \ncm{\OExt}{\mbox{\rm OrderExt}}
\ncm{\AI}{\mbox{\rm AInn($A$)}} \ncm{\HI}{\mbox{\rm HInn($A$)}}
\ncm{\Aut}{\mbox{\rm Aut}} \ncm{\Mal}{\mbox{$M_{\alpha}$}}
\ncm{\Aff}{\mbox{${\rm Aff}$}} \ncm{\id}{\mbox{\rm id}}
\ncm{\Ker}{\mbox{\rm Ker}} \ncm{\BE}{\begin{eqnarray*}}
\ncm{\EE}{\end{eqnarray*}} \ncm{\lra}{\mbox{$\longrightarrow$}}
\ncm{\Hom}{\mbox{\rm Hom}} \ncm{\calU}{{\cal U}} \ncm{\el}{\ell}
\ncm{\ad}{\mbox{\rm ad}} \ncm{\Alg}{\mbox{\rm Alg}}
\ncm{\Conv}{\mbox{\rm Conv}}\ncm{\D}{{\cal D}}

\ncm{\cstar}{$C^{*}$-algebra} \ncm{\cstars}{$C^{*}$-algebras}
\ncm{\ra}{\mbox{$\rightarrow$}} \ncm{\la}{\mbox{$\leftarrow$}}
\ncm{\hra}{\hookrightarrow} \ncm{\da}{\mbox{$\downarrow$}}
\ncm{\se}{\mbox{$\searrow$}} \ncm{\al}{\mbox{$\alpha $}}
\ncm{\del}{\mbox{$\delta$}} \ncm{\supp}{\mbox{\rm supp}}
\ncm{\Ad}{\mbox{\rm Ad}} \ncm{\CAR}{\mbox{$M_{2^{\infty}}$}}
\ncm{\ep}{\mbox{$\epsilon > 0$}} \ncm{\mod}{\mbox{\rm mod}}
\ncm{\Sp}{\mbox{\rm Sp}} \ncm{\ol}{\overline}
\ncm{\Mninf}{\mbox{$M_{n^{\infty}}$}} \ncm{\MR}{M. R\o{}rdam}
\ncm{\Range}{\mbox{\rm Range}}
\ncm{\vo}{}
\ncm{\ch}{}
\ncm{\CMP}{Comm. Math. Phys.} \ncm{\add}{}
\ncm{\tilsig}{\tilde{\sigma}}

\newtheorem{theo}{Theorem}[section]

\newtheorem{lem}[theo]{Lemma}
\newtheorem{prop}[theo]{Proposition}
\newtheorem{remark}[theo]{Remark}
\newtheorem{definition}[theo]{Definition}
\newtheorem{example}[theo]{Example}

\newenvironment{rem}{\begin{remark} \rm}{\end{remark}}
\newenvironment{pf}{{\it Proof.}}{\vspace{3mm}}

\ncm{\R}{\mbox{\bf R}} \ncm{\Z}{\mbox{\bf Z}} \ncm{\T}{\mbox{\bf
T}} \ncm{\TT}{\T$^{2}$} \ncm{\N}{\mbox{\bf N}} \ncm{\C}{\mbox{\bf
C}}

\title{AF flows and continuous symmetries}

\bigskip

\author{O. Bratteli\\
   {\small Department of Mathematics,
   University of Oslo, PB 1053--Blindern, N-0316 Oslo, Norway}\\
                A. Kishimoto\\
   \small Department of Mathematics, Hokkaido University,
    Sapporo 060, Japan}
\begin{document}
\maketitle

\begin{abstract}
We consider AF flows, i.e., one-parameter automorphism groups of a
unital simple AF $C^*$-algebra which leave invariant the dense
union of an increasing sequence of finite-dimensional
*-subalgebras, and derive two properties for these; an absence of
continuous symmetry breaking and a kind of  real rank zero
property for the {\em almost fixed points}.
\end{abstract}

\section{Introduction}
We consider the class of AF representable one-parameter
automorphism groups of a unital simple AF \cstar\ (which will be
called AF flows) and derive two properties, one of which is
invariant under inner perturbations and may be used to distinguish
them from other flows (i.e., one-parameter automorphism groups).

We recall that a flow $\alpha$ of a unital simple AF \cstar\ $A$
is defined to be AF locally representable or an AF flow if there is
an increasing sequence $(A_n)$ of $\alpha$-invariant
finite-dimensional *-subalgebras of $A$ with dense union
\cite{Kis99a}, \cite{Kis99}. In this case there is a self-adjoint
$h_n\in A_n$ such that
$\alpha_t|A_n=\Ad\,e^{ith_n}|A_n$ for each $n$. Thus the {\em
local Hamiltonians} $(h_n)$ mutually commute and can be considered
to represent the time evolution of a classical statistical lattice
model, which is a special kind of model among all the models
quantum or classical. Consider the larger class of flows which are
inner perturbations of AF-flows. (These are characterized by the
property that the domains of the generators contains a canonical AF
maximal abelian sub-algebra (masa), see
\cite[Proposition~3.1]{Kis99}.) In \cite[Theorem~2.1 and
Remark~3.3]{Kis99} it was demonstrated that there are flows
outside this larger class, but the proof was not easy.
Our original aim was to
show that all the flows which naturally arose in quantum
statistical lattice models and were not obviously AF flows, were
in fact beyond the class of inner perturbations of AF flows. We
could not prove that there was even a single example and obtained
only a weak result in this direction which is presented in
Remark~2.4. Thus we ended up presenting the two new properties of
the AF flows mentioned in the abstract.

The first property we derive for AF flows can be expressed as:
{\em there is no continuous symmetry breaking}. If
$\delta_{\alpha}$ denotes the generator of a general flow $\alpha$,
we define the {\em exact} symmetry group for $\alpha$ as
$G_0=\{\gamma\in\Aut A\ |\
\gamma\delta_{\alpha}\gamma^{-1}=\delta_{\alpha}\}$ and the {\em
near} symmetry group as $G_1=\{\gamma\in\Aut A\ |\
\gamma\delta_{\alpha}\gamma^{-1}=\delta_{\alpha}+\ad\,ih\ {\rm
for\ some}\ h=h^*\in A\}$. Then it is known that there is a
natural homomorphism of $G_0$ into the affine homeomorphism group
of the simplex of KMS states at each temperature. We deduce
moreover in Proposition~2.1 from the perturbation theory of KMS
states
\cite{Ar}, that there is a homomorphism of $G_1$ into the
homeomorphism group of the simplex of KMS states at each
temperature, mapping the extreme points onto the extreme points. We
next show in the special case of AF flows that if $\gamma\in G_0$ is
connected to
$\id$ in
$G_0$ by a continuous path, then
$\gamma$ induces the identity map on the simplexes of KMS states.
We actually show a generalization of this in Theorem~2.3: If
$\alpha$ is an AF flow and $\gamma\in G_1$ is connected to $\id$ in
$G_1$ by a continuous path
$(\gamma_t)$ such that
$\gamma_t\delta_{\alpha}\gamma_t^{-1}=\delta_{\alpha}+\ad\,ib(t)$
with $b(t)$ rectifiable in $A$, then $\gamma$ induces a
homeomorphism which fixes each extreme point. (Thus, if the
homeomorphism is affine, it is the identity map. This is in
particular true if $\gamma\in G_0$.)

The second property we derive for the class of inner perturbations
of AF flows can be expressed as: {\em the almost fixed point
algebra for $\alpha$ has real rank zero} (see Theorem \ref{C4}). A
technical lemma used to show this property is a generalization of
H. Lin's result on almost commuting self-adjoint matrices
\cite{L}. The generalization says that any almost commuting pair
of self-adjoint matrices, one of norm one and the other of
arbitrary norm, is in fact close to an exactly commuting pair (see
Theorem~\ref{C}).

We recall here a similar kind of property in \cite{Kis99} saying
that {\em the almost fixed point algebra has trivial $K_1$}. We
will show by examples that these two properties, real rank zero
and trivial $K_1$ for the almost fixed point algebra, are
independent, as one would expect. (It is not that {\em the almost
fixed point algebra} is actually defined as an algebra; but if
$\alpha$ is periodic, then we can regard the almost fixed point
algebra as the usual fixed point algebra, see Proposition~3.7. In
general we can characterize any property of {\em the almost fixed
point algebra\/} as the corresponding property of the fixed point
algebra for a certain flow obtained by passing to a \cstar\ of
bounded sequences modulo
$c_0$, see Proposition~3.8.)

We remark that there is a flow $\alpha$ of a unital simple AF
$C^*$-algebra such that $\D(\delta_{\alpha})$ is not AF (as a
Banach *-algebra)(cf. \cite{Sak0,Sak}). This was shown in
\cite{Kis99} by constructing an example where
$\D(\delta_{\alpha})$ does not have real rank zero. Note that
$\D(\delta_{\alpha})$ has always trivial $K_1$ and has the same
$K_0$ as the \cstar\ $A$. Hence real rank is still the only
property which has been used to distinguish $\alpha$ with non-AF
$\D(\delta_{\alpha})$. On the other hand even $K_0$ (of the
almost fixed point algebra) might be used  to distinguish non-AF
flows (up to inner perturbations) as well as real rank and $K_1$
as shown above.

In the last section we will show that any quasi-free flow of the
CAR algebra has the property that the almost fixed point
algebra has trivial $K_1$, leaving open the question of whether
it is an inner perturbation of an AF flow or not and even the weaker
question of whether the almost fixed point algebra has real rank
zero or not.

One of the authors (A.K.) would like to thank Professor S.~Sakai
for discussions and questions concerning the first property.

\section{Symmetry}

In the first part of this section we describe the symmetry group
of a flow and how it is mapped into the homeomorphism groups of
the simplexes of KMS states. Then in the remaining part we discuss
a theorem on a kind of absence of continuous symmetry
breaking for AF flows.

In the first part $A$ can be an arbitrary unital simple \cstar.
Let $\alpha$ be a flow of $A$ (where we always assume strong
continuity; $t\mapsto \alpha_t(x)$ is continuous for any $x\in
A$), and $\delta_{\alpha}$ the generator of $\alpha$. Then
$\delta_{\alpha}$ is a closed linear operator defined on a dense
*-subalgebra $\D(\delta_{\alpha})$ of $A$ with the derivation
property:
$\delta_{\alpha}(xy)=\delta_{\alpha}(x)y+x\delta_{\alpha}(y),\
\delta_{\alpha}(x)^*=\delta_{\alpha}(x^*)$ for $x,y\in
\D(\delta_{\alpha})$. We equip $\D(\delta_{\alpha})$ with the norm
$\|\,\cdot\,\|_{\delta_{\alpha}}$ obtained by embedding
$\D(\delta_{\alpha})$ into $A\otimes M_2$ by the (non *-preserving)
isomorphism $x\mapsto
\left(\begin{array}{cc}x&\delta_{\alpha}(x)\\0&x\end{array}\right)$.
Note that $\D(\delta_{\alpha})$ is a Banach *-algebra.  (See
\cite{BR1,Br,Sak} for the theory of unbounded derivations.)

We call a continuous function $u$ of $\R$ into the unitary group
of $A$ an $\alpha$-cocycle if $u_s\alpha_s(u_t)=u_{s+t},\
s,t\in\R$. Then $t\mapsto \Ad\,u_t\circ\alpha_t$ is a flow of $A$
and is called a cocycle perturbation of $\alpha$. If $u$ is
differentiable, then the generator of this perturbation is
$\delta_{\alpha}+\ad\,ih$, where $du/dt|_{t=0}=ih$ (see
\cite[section~1]{Kis99a}). We define the symmetry group
$G=G_{\alpha}$ of
$\alpha$  as
$$
\{\gamma\in\Aut\,A\ |\ \gamma\alpha\gamma^{-1}\ {\rm is\ a\
cocycle\ perturbation\ of}\ \alpha\},
$$
which is slightly more general than the $G_1$ given in the
introduction, so $G_0\subseteq G_1\subseteq G=G_\alpha$. Then $G$
depends on the class of cocycle perturbations of $\alpha$ only and
is indeed a group: If
$\gamma\in G$, then $\gamma\alpha_t\gamma^{-1}=\Ad\,u_t\,\alpha_t$
for some $\alpha$-cocycle $u$, which implies that
  $$
  \gamma^{-1}\alpha_t\gamma=\Ad\gamma^{-1}(u_t^*)\alpha_t.
  $$
We can check the $\alpha$-cocycle property of $t\mapsto
\gamma^{-1}(u_t^*)$ by
\BE
\gamma^{-1}(u_s^*)\alpha_s(\gamma^{-1}(u_t^*))&=&
\gamma^{-1}(u_s^*\gamma\alpha_s\gamma^{-1}(u_t^*))
=\gamma^{-1}(u_s^*\Ad\,u_s\,\alpha_s(u_t^*))\\
&=&\gamma^{-1}(\alpha_s(u_t^*)u_s^*)=\gamma^{-1}(u_{s+t}^*).
\EE
Thus $\gamma^{-1}\in G$. If $\gamma_1,\ \gamma_2\in G$, then
$\gamma_i\alpha_t\gamma_i^{-1}=\Ad\,u_{it}\,\alpha_t$ for some
$\alpha$-cocycle $u_i$ for $i=1,2$. Since
$\gamma_1\gamma_2\alpha_t(\gamma_1\gamma_2)^{-1}=\Ad\,\gamma_1(u_{2t})u_{1t}
\,\alpha_t$,
we only have to check that $t\mapsto \gamma_1(u_{2t})u_{1t}$ is an
$\alpha$-cocycle, which will be denoted by $\gamma_1(u_2)u_1$. We
leave this simple calculation to the reader. Note that $G$
contains the inner automorphism group ${\rm Inn}(A)$ as a normal
subgroup and each element of $G/{\rm Inn}(A)$ has a representative
$\gamma\in G$ such that $\gamma$ leaves $\D(\delta_{\alpha})$
invariant and
$$
\gamma\delta_{\alpha}\gamma^{-1}=\delta_{\alpha}+\ad\,ib
$$
for some $b=b^*\in A$ (see \cite[Corollary~1.2]{Kis99a}).

We equip $G=G_{\alpha}$ with the topology defined by
$\gamma_n\ra\gamma$ in $G$ if
\begin{itemize}
\item[(1)]
$\gamma_n\ra \gamma$ in $\Aut(A)$
(i.e., $\|\gamma_n(x)-\gamma(x)\|\ra0$ for $x\in A$),
\end{itemize}
and
\begin{itemize}
\item[(2)]
there
exist $\alpha$-cocycles $u_n,\ u$ such that
$\gamma_n\alpha_t\gamma_n^{-1}=\Ad\,u_{nt}\,\alpha_t$,
$\gamma\alpha_t\gamma^{-1}=\Ad u_t\alpha_t$ and
$\|u_{nt}-u_t\|\ra0$ uniformly in $t$ on compact subsets of $\R$.
\end{itemize}

\noindent
With this topology $G$ is a topological group.

Let $c\in\R\setminus \{0\}$ and $\omega$ a state on $A$. We say
that $\omega$ satisfies the $c$-KMS condition or is a $c$-KMS
state (with respect to $\alpha$) if for any $x,y\in A$ there is a
bounded continuous function $F$ on the strip $S_c=\{z\in\C\ |\
0\leq \Im z/c\leq 1\}$ such that $F$ is analytic in the interior
of $S_c$ and satisfies, on the boundary of $S_c$,
  \BE F(t)&=& \omega(x\alpha_t(y)),\ \ t\in\R,\\
  F(t+ic)&=&\omega(\alpha_t(y)x),\ \ t\in \R.
  \EE
We denote by $K_c^{\alpha}=K_c$ the set of $c$-KMS states of $A$.
Then $K_c$ is a closed convex set of states and moreover a
simplex. We denote by $\partial(K_c)$ the set of extreme points of
$K_c$. Note that for $\omega\in K_c$, $\omega$ is extreme in $K_c$
if and only if $\omega$ is a factorial state (see \cite{BR1,Sak}
for details).

\begin{prop}\label{A}
Let $A$ be a unital simple \cstar, $\alpha$ a flow of $A$, and
$c\in\R\setminus\{0\}$. Then there is a continuous homomorphism
$\Phi$ of the symmetry group $G_{\alpha}$ of $\alpha$ into the
homeomorphism group of $K_c$ such that $\Phi(\gamma)(\omega)$ is
unitarily equivalent to $\omega\gamma^{-1}$ for each $\gamma\in
G_{\alpha}$ and $\omega\in K_c$. Moreover
$\Phi(\gamma)=\id$ for any inner $\gamma$.
\end{prop}
\begin{pf}
Let $\gamma\in G_{\alpha}$ and let $u$ be an $\alpha$-cocycle such
that $\gamma\alpha_t\gamma^{-1}=\Ad\,u_t\,\alpha_t$. Since $A$ is
simple, $u$ is unique up to phase factors, i.e., any other
$\alpha$-cocycle satisfying the same equality is given as
$t\mapsto e^{ipt}u_t$ for some $p\in\R$.

Let $\omega\in K_c$. Then $\omega\gamma^{-1}$ is a KMS state with
respect to
$\gamma\alpha_t\gamma^{-1}=\Ad\,u_t\,\alpha_t$. Using the fact
that $\alpha_t=\Ad\,u_t^*\,\gamma\alpha_t\gamma^{-1}$, there is a
procedure to make a KMS positive linear functional $\omega'$ with
respect to $\alpha$, which depends on the choice of $u$; formally
it can be given as
  $$
  \omega'(x)=\omega\gamma^{-1}(xu_{ic}^*),\ \ x\in A.
  $$
More precisely we let $\beta_t=\Ad\,u_t\,\alpha_t$ and express the
$\beta$-cocycle $u_t^*$ as
  $$
  u_t^*=wv_t\beta_t(w^{-1})
  $$
such that $t\mapsto v_t$ extends to an entire function on $\C$
\cite[Lemma~1.1]{Kis99a}. Then we define $\Phi(\gamma,u)\omega$ as
  $$
  (\Phi(\gamma,u)(\omega))(a)=\omega\gamma^{-1}(w^{-1}awv_{ic}).
  $$
(By a formal calculation we can see that this satisfies the
$c$-KMS condition as follows:
  \BE
  \omega\gamma^{-1}(w^{-1}a\alpha_{ic}(b)wv_{ic})
  &=&\omega\gamma^{-1}(w^{-1}au_{ic}^{-1}
  \beta_{ic}(b)u_{ic}wv_{ic})\\
  &=&\omega\gamma^{-1}(w^{-1}awv_{ic}\beta_{ic}(w^{-1}bw))\\
  &=&\omega\gamma^{-1}(w^{-1}bawv_{ic}),
  \EE
where we used that $u_{ic}=\beta_{ic}(w)v_{ic}^{-1}w^{-1}$ and that
$\omega\gamma^{-1}$ is a $c$-KMS state for
$\beta_t=\gamma\alpha_t\gamma^{-1}$. See
\cite{Kis99a}.)  The map
$$
\Phi(\gamma):\omega\mapsto
\Phi(\gamma,u)(\omega)/\Phi(\gamma,u)(\omega)(1)
$$
defines a continuous map of $K_c$ into $K_c$ and
$\Phi(\gamma,u)(\omega)$ is quasi-equivalent (hence unitarily
equivalent) to $\omega\gamma^{-1}$.
(It follows from the definition of $\Phi(\gamma,v)$ that
$\Phi(\gamma)(\omega)$ is quasi-contained in $\omega\gamma^{-1}$,
but as $w^{-1}$ and $wv_{ic}$ are invertible,
$\omega\gamma^{-1}$ is conversely quasi-contained in
$\Phi(\gamma)(\omega)$. Since any KMS state is separating and
cyclic for the weak closure, these states are unitary equivalent.)
 For any other choice
$u_t'=e^{ipt}u_t$ for $u$ it follows that
$\Phi(\gamma,u')=e^{-cp}\Phi(\gamma,u)$. Thus $\Phi(\gamma)$ does
not depend on the choice of $u$. For $\gamma_1,\gamma_2\in
G_{\alpha}$ with $\alpha$-cocycles $u_1,u_2$ respectively, it
follows that
$$
\Phi(\gamma_1\gamma_2,\gamma_1(u_2)u_1)=
\Phi(\gamma_1,u_1)\Phi(\gamma_2,u_2)
  $$
since
  \BE\Phi(\gamma_1,u_1)\Phi(\gamma_2,u_2)(\omega)
  &=&\Phi(\gamma_1,u_1)(\omega\gamma_2^{-1}(\,\cdot\,u_{2,ic}^{*}))\\
  &=&
    \omega\gamma_2^{-1}(\gamma_1^{-1}(\,\cdot\,u_{1,ic}^*)
    u_{2,ic}^{*})\\
&=&\omega\gamma_2^{-1}\gamma_1^{-1}(\,\cdot\,u_{1,ic}^*
\gamma_1(u_{2,ic}^*)).
  \EE
This shows that $\Phi$ is a group homomorphism. If
$\gamma=\Ad\,u$, then $\Phi(\gamma,u\alpha(u^*))(\omega)=\omega$.
The continuity of $\gamma\mapsto \Phi(\gamma)$ follows from the
following lemma.
\end{pf}

\begin{lem}
Let $(u_{\infty},u_1,u_2,\ldots)$ be a sequence of
$\alpha$-cocycles such that
$\lim_{n\rightarrow\infty}u_{n,t}=u_{\infty,t}$ uniformly in $t$
on every compact subset of \R. Then for any $\epsilon>0$ there
exists a sequence $(w_{\infty},w_1,w_2,\ldots)$ of invertible
elements in $A$ such that
$\lim_{n\rightarrow\infty}w_n=w_{\infty},\ \|w_n-1\|<\epsilon$,
and $v_{m,t}\equiv w_{m}^{-1}u_{m,t}\alpha_t(w_m)$ extends to an
entire function on $\C$ for $m=\infty,1,2\ldots$ such that
$\lim_{n\rightarrow\infty}v_{n,z}=v_{\infty,z}$ for any $z\in\C$.
\end{lem}
\begin{pf}
Define a \cstar\ $B$ by
  $$
  B=\{x=(x_n)_{n=1}^{\infty}\ |\ x_n\in A,\ \lim x_n\
\mbox{exists}\}
  $$
and define a flow $\beta$ on $B\otimes M_2$ by $
\beta_t=\Ad\,U\circ \alpha_t\otimes\id$, where $U=(1\oplus
u_{n,t})$. We define a homomorphism $\varphi$ of $B$ onto $A$ by
$\varphi(x)=\lim x_n$ for $x=(x_n)\in B$ and note that
$\varphi\circ\beta_t=\Ad(1\oplus
u_{\infty,t})\circ\alpha_t\otimes\id\circ\varphi$. Let
$\epsilon\in (0,1)$. Since $(1\oplus0)_n$ and $(0\oplus 1)_n$ are
fixed by $\beta$, there is a $w\in B$ such that $\|w-1\|<\epsilon$
and
  $$
  t\mapsto \beta_t(\left(\begin{array}{cc}0&0\\w&0\end{array}\right))
  $$
extends to an entire function on $\C$ (pick an entire element $y$
for $\beta$ close to
$\left(
{0\;\,0 \atop 1\;\,0}\right)$,
and replace
$y$ by
$(0\otimes 1)_n y(1\otimes 0)_n$). If
$w=(w_n)\in B,\ v_{n,t}= w_n^{-1}u_{n,t}\alpha_t(w_n)\in A$, and
$v_t=(v_{n,t})\in B$, then we have that
  $$
  \beta_t(\left(\begin{array}{cc}0&0\\w&0\end{array}\right))=
  \left(\begin{array}{cc}0&0\\wv_t&0\end{array}\right).
  $$
Letting $w_{\infty}=\lim w_n$ and $v_{\infty,t}=\lim v_{n,t}$,
the proof is complete.
\end{pf}

\begin{theo}\label{B}
Let $A$ be a unital simple AF \cstar\ and $\alpha$ an AF flow of
$A$. Let $(\gamma_t)_{t\in [0,1]}$ be a continuous path in
$G_{\alpha}$ such that
  $$
  \gamma_t\delta_{\alpha}\gamma_t^{-1}=\delta_{\alpha}+\ad\,ib(t)
  $$
for some rectifiable path $(b(t))_{t\in [0,1]}$ in $A_{sa}$. Then
it follows that $\Phi(\gamma_0)(\omega)=\Phi(\gamma_1)(\omega)$
for $\omega\in\partial(K_c)$.
\end{theo}
\begin{pf}
Let $C$ be a canonical AF masa in $\D(\delta_{\alpha})$ such that
$\delta_{\alpha}|_C=0$. Let $\omega\in K_c$. We note that if $E$
denotes the projection of norm one onto $C$, then
$\omega=(\omega\big|_C)\circ E$, i.e., $\omega$ is determined by the
restriction $\omega\big|_C$. (Let $(A_n)$ be an increasing sequence
of
$\alpha$-invariant finite dimensional subalgebras with dense union
in $A$ such that $A_n\cap C$ is masa in $A_n$ for each $n$. Then
$\omega\big|_{A_n}$ is clearly determined by $\omega\big|_{A_n\cap
C}$, and thus
$\omega$ is determined by $\omega\big|_C$.)

We first prove the theorem in the simpler case where $b(t)=0$. In
this case $\gamma_t$ leaves the $C^*$-subalgebra $B={\rm
Kernel}(\delta_{\alpha})$ invariant, on which $\omega$ is a trace.
For any projection $e\in C\subset B$, $(\gamma_t(e))$ is a
continuous family of projections in $B$, which implies that
$\gamma_0(e)$ is equivalent to $\gamma_1(e)$ in $B$. Hence
$\omega\gamma_0(e)=\omega\gamma_1(e)$. Since $C$ is an abelian AF
algebra, this implies that $\omega\gamma_0|_C=\omega\gamma_1|_C$.
Since they are KMS states, we can conclude that
$\omega\gamma_0=\omega\gamma_1$. Since this is true for any
$\omega\in K_c$, it also follows that
$\omega\gamma_0^{-1}=\omega\gamma_1^{-1}$.

What we will do in the following is a modification of this
argument.

Let $\omega\in \partial(K_c)$. In the GNS representation associated
with $\omega\in \partial(K_c)$, we define a one-parameter unitary
group
$U$ by
  $$
U_t\pi_{\omega}(x)\Omega_{\omega}=\pi_{\omega}\circ
    \alpha_t(x)\Omega_{\omega},
  \ \ \ x\in A.
  $$
Then from the $c$-KMS condition on $\omega$ it follows that the
modular operator $\Delta$ for $\Omega_{\omega}$ is given by
$\Delta=e^{-cH}$, where $H$ is the generator of $U$;
$U_t=e^{itH}$ (See \cite[Proof of Theorem~5.3.10]{BR2}). We define a
positive linear functional
$\omega^{(h)}$ on $A$ for $h=h^*\in A$ as the vector state given
by $e^{-c(H+\pi_{\omega}(h))/2}\Omega_{\omega}$, i.e.,
  $$
  \omega^{(h)}(x)=(\pi_{\omega}(x)e^{-c(H+\pi_{\omega}(h))/2}
   \Omega_{\omega},
  e^{-c(H+\pi_{\omega}(h))/2}\Omega_{\omega}).
  $$
Then $\omega^{(h)}$ satisfies the $c$-KMS condition with respect
to $\delta_{\alpha}+\ad\,ih$. (See \cite{Ar,Sak} or
\cite[Theorem~5.4.4]{BR2}. The relation to the previous
perturbation argument in terms of cocycles is as follows: The flow
generated by
$\delta_{\alpha}+\ad\,ih$ is given as $\Ad\,u_t\,\alpha_t$, where
$u$ is the $\alpha$-cocycle with
$du_t/dt|_{t=0}=ih$, and $\omega^{(h)}$ is equal to
$\omega(w^{-1}\,\cdot\,wv_{ic})$, where $u_t$ is expressed as
$wv_t\alpha_t(w^{-1})$ with $t\mapsto v_t$ entire.)

For $s\in [0,1]$ let $\omega_s=\omega^{(b(s))}$, which is a
positive linear functional satisfying the $c$-KMS
condition with respect to the generator
$\delta_{\alpha}+\ad\,ib(s)$. This implies that $\omega_s\gamma_s$
is a $c$-KMS positive linear functional with respect to
$\gamma_s^{-1}(\delta_{\alpha}+\ad\,ib(s))\gamma_s=\delta_{\alpha}$.

Let $s_1,s_2\in [0,1]$ and define a positive linear functional
$\varphi$ on
$A\otimes M_2$ by
  $$
  \varphi(a)=\omega_{s_1}(a_{11})+\omega_{s_2}(a_{22})
  $$
for $a=(a_{ij})\in A\otimes M_2$. Then $\varphi$ is a $c$-KMS
positive linear functional for the flow $\beta$ of $A\otimes M_2$
defined by
  $$
  \beta_t((a_{ij}))=\Ad\left(\begin{array}{cc}u_t^{(b(s_1))}&0\\
                                   0&u_t^{(b(s_2))}\end{array}\right)
                    (\alpha_t(a_{ij})),
  $$
where $u_t^{(h)}$ is the $\alpha$-cocycle determined by
$du_t^{(h)}/dt|_{t=0}=ih$ (see \cite{Connes}). The generator
$\delta_{\beta}$ of $\beta$ is given by
  $$
  \delta_{\beta}((a_{ij}))
  =\left(\begin{array}{cc}(\delta_{\alpha}+\ad\,ib(s_1))(a_{11})&
       \delta_{\alpha}(a_{12})+ib(s_1)a_{12}-a_{12}ib(s_2)\\
      \delta_{\alpha}(a_{21})-a_{21}ib(s_{1})+ib(s_2)a_{21}&
     (\delta_{\alpha}+\ad\,ib(s_2))(a_{22})\end{array}\right).
  $$

Fix $\epsilon\in (0,1/2)$ and a $C^{\infty}$-function $f$ on $\R$
with compact support such that $f(0)=0$ and $f(t)=t^{-1/2}$ on
$[1-\epsilon,1]$. Let $e$ be a projection in $C$. We choose
$s_1,s_2\in [0,1]$ so that
  $$
  \|\gamma_{s_1}(e)-\gamma_{s_2}(e)\|<\epsilon.
  $$
Let
  $$
  x=\left(\begin{array}{cc}0&\gamma_{s_1}(e)\gamma_{s_2}(e)\\
                     0&0\end{array}\right).
  $$
Then
  $$
  x^*x=\left(\begin{array}{cc}
  0&0\\0&\gamma_{s_2}(e)\gamma_{s_1}(e)\gamma_{s_2}(e)\end{array}
  \right)
  $$
and $\Sp(x^*x)\subset \{0\}\cup (1-\epsilon,1]$. Let $v=xf(x^*x)$.
Then $v$ is a partial isometry such that
  $$
  vv^*=\left(\begin{array}{cc}\gamma_{s_1}(e)&0\\0&0\end{array}
   \right),\ \ \
  v^*v=\left(\begin{array}{cc}0&0\\0&\gamma_{s_2}(e)\end{array}
   \right).
  $$
Since all the components of $\delta_{\beta}(x)$ are zero except
for the (1,2) component and
$(\delta_{\alpha}+\ad\,ib(s))\gamma_s(e)=0$, we have that
  \BE
  \|\delta_{\beta}(x)\|&=& \|\delta_{\beta}(x)_{12}\|\\
   &=&
\|\gamma_{s_1}(e)ib(s_1)\gamma_{s_2}(e)-\gamma_{s_1}ib(s_2)
  \gamma_{s_2}(e)\| \\
   &\leq& \|b(s_1)-b(s_2)\|.
  \EE
Since $\|\delta_{\beta}(x^*x)\|\leq 2\|b(s_1)-b(s_2)\|$, and
  $$
  \delta_{\beta}(f(x^*x))=\delta_{\beta}(\int\hat{f}(s)
    e^{isx^*x}ds)=
  \int\hat{f}(s)\int_0^1e^{itsx^*x}is\delta_{\beta}
    (x^*x)e^{i(1-t)sx^*x}dtds,
  $$
it follows that
  $$
  \|\delta_{\beta}(f(x^*x))\|\leq
   \int|\hat{f}(s)s|ds\cdot \|\delta_{\beta}(x^*x)\|.
  $$
Thus there is a constant $C>0$ such that
  $$
  \|\delta_{\beta}(v)\|\leq C\|b(s_1)-b(s_2)\|.
  $$

By the KMS condition on $\varphi$ we have a continuous function
$f$ on the strip $S_c$ between $\Im z=0$ and $\Im z=c$, analytic
in the interior, such that
  \BE
  f(t)&=& \varphi(v\beta_t(v^*)),\ \ \ t\in \R,\\
  f(t+ic)&=&\varphi(\beta_t(v^*)v),\ \ \ t\in\R.
  \EE
Then $f$ is differentiable on $S_c$ including the boundary and
satisfies that
  \BE
  f'(t)&=&\varphi(v\beta_t(\delta_{\beta}(v^*))),\ \ \ t\in\R,\\
  f'(t+ic)&=&\varphi(\beta_t(\delta_{\beta}(v^*))v), \ \ \  t\in \R.
  \EE
Hence it follows that
  $$
  \sup_{z\in S_c}|f'(z)|\leq \sup_{z\in\partial S_c}|f'(z)|\leq C
  \max\{\|\omega_{s_1}\|,\|\omega_{s_2}\|\}\|b(s_1)-b(s_2)\|,
  $$
which implies that
  \BE
  |\omega_{s_2}(\gamma_{s_2}(e))-\omega_{s_1}(\gamma_{s_1}(e))|&=&
   |f(ic)-f(0)|\\
  &\leq& |c|CM\|b(s_1)-b(s_2)\|,
  \EE
where $M=\max\{\|\omega_s\|\ |\ s\in [0,1]\}$. We let
$m=\min\{\|\omega_s\|\ |\ s\in [0,1]\}$ and choose
$t_0=0<t_1<\cdots t_k=1$ such that
  $$
  |c|C\frac{M}{m}(1+\frac{M}{m}){\rm Length}(b(s),\ s\in[t_{i-1},t_i])
  <\frac{1}{4}.
  $$
Then for any projection $e\in C$, we subdivide each interval
$[t_{i-1},t_i]$ into $s_0=t_{i-1}<s_1<\cdots <s_{\ell}=t_i$ such
that
  $$
  \|
  \gamma_{s_{j-1}}(e)-\gamma_{s_j}(e)\|<\epsilon,
  $$
and apply the above argument to each pair $s_{j-1},s_j$ to obtain
that
  $$
  |\omega_{t_{i-1}}\gamma_{t_{i-1}}(e)-\omega_{t_i}\gamma_{t_i}(e)|\leq
  |c|CM\,{\rm Length}(b(s),\ s\in [t_{i-1},t_i]).
  $$
Thus we have that for any projection $e\in C$
  \BE
  |\frac{\omega_{t_{i-1}}\gamma_{t_{i-1}}(e)}{\omega_{t_{i-1}}(1)}-
   \frac{\omega_{t_i}\gamma_{t_i}(e)}{\omega_{t_i}(1)}|
   &\leq&
    \frac{1}{m}|\omega_{t_{i-1}}
    \gamma_{t_{i-1}}(e)-\omega_{t_i}\gamma_{t_i}(e)|+
   M|\frac{\omega_{t_i}(1)-\omega_{t_{i-1}}(1)}
   {\omega_{t_i}(1)\omega_{t_{i-1}}(1)}|\\
   &\leq&(\frac{1}{m}+\frac{M}{m^2})
   |c|CM{\rm Length}(b(s),\ s\in [t_{i-1},t_i])\\
   &\leq& 1/4.
  \EE
Let
  $$
  \varphi_t=\frac{\omega_t\gamma_t}{\omega_t(1)}
  $$
and recall that $\varphi_t$ is a factorial $c$-KMS state with
respect to $\alpha$. Since $\varphi_t=\varphi_t E$ with $E$ the
projection onto $C$ and
$\|(\varphi_{t_{i-1}}-\varphi_{t_i})|C\|\leq1/2$, we have that
$\|\varphi_{t_{i-1}}-\varphi_{t_i}\|\leq1/2$. Hence
$\varphi_{t_{i-1}}=\varphi_{t_i}$. Thus we conclude that
$\varphi_0=\varphi_1$ or $
  \Phi(\gamma_0^{-1})(\omega)=\Phi(\gamma_1^{-1})(\omega)
  $
for $\omega\in\partial(K_c)$. This implies that
$\Phi(\gamma_0)(\omega)=\Phi(\gamma_1)(\omega)$ for
$\omega\in\partial(K_c)$ as well.
\end{pf}

\begin{rem}
Among the quantum lattice models, two or more dimensional, there
are long-range interactions which exhibit continuous symmetry
breaking. Let $\alpha$ be the flow generated by such an interaction
and let $\gamma$ be an action of $\T$ which exactly commutes with
$\alpha$ and acts non-trivially on the simplex of $c$-KMS states
at some inverse temperature $c>0$. Suppose that $\alpha$ is an
inner perturbation of an AF flow, i.e.,
$\delta=\delta_{\alpha}+\ad\,ib$ is the generator of an AF flow.
Since $\gamma_t\delta\gamma_t^{-1}=\delta+\ad\,i(\gamma_t(b)-b)$,
we can conclude that $t\mapsto \gamma_t(b)$ is not rectifiable;
thus at least $b$ is not in the domain of the generator of
$\gamma$. (Note we still cannot conclude that $\alpha$ is not an
inner perturbation of an AF flow.) \end{rem}

\section{Property of real rank zero}
\setcounter{theo}{0}

First we generalize H. Lin's result \cite{L}
and then use it to prove that {\em the almost fixed point
algebra for an AF flow has real rank zero}.

\begin{theo}\label{C}
For every $\epsilon>0$ there is a $\nu>0$ satisfying the following
condition: For any $n\in\N$ and any pair $a,b\in (M_n)_{sa}$ with
$\|b\|\leq 1$ and $\|[a,b]\|<\nu$ there exists a pair $a_1,b_1\in
(M_n)_{sa}$ such that $\|a-a_1\|<\epsilon,\ \|b-b_1\|<\epsilon$,
and $[a_1,b_1]=0$.
\end{theo}

If we impose the extra condition that $\|a\|\leq1$ for $a$, then
this result is due to H. Lin (see also \cite{FR}). Our proof is to
reduce Theorem~3.1 to Lin's result.

\begin{lem}\label{C1}
Let $f$ be a $C^{\infty}$-function on $\R$ such that $f\geq0,\
\int f(t)dt=1$, and $\supp\,\hat{f}\subset (-1/2,1/2)$. For any
pair
$a,b$ elements in a \cstars-algebra such that $a=a^\ast$, define
  $$
  b_1=\int f(t)e^{ita}be^{-ita}dt.
  $$
Then it follows that
  \BE
  \|b-b_1\|&\leq& \int f(t)|t|dt\cdot \|[a,b]\|,\\
  \|[a,b_1]\|&\leq& \int f(t) d t\cdot \|[a,b]\|.
  \EE
\end{lem}
\begin{pf}
This follows from the following computations:
  \BE
  b_1-b&=&\int f(t)(e^{ita}be^{-ita}-b)dt,\\
   &=&\int f(t)\int_0^t  e^{isa}[ia,b]e^{-isa}dsdt,\\
   [0pt]
    [ a, b_1]&=&\int f(t)e^{ita}[a,b]e^{-ita}dt.
  \EE
\end{pf}

\begin{rem}
If we denote by $E_a$ the spectral measure of $a$, then the $b_1$
defined in the above lemma satisfies that
  $$
  E_a(-\infty,t-1/4]\,b_1E_a[t+1/4,\infty)=0
  $$
for any $t\in\R$, \cite[Proposition~3.2.43]{BR1}.
\end{rem}

\begin{lem}\label{C2}
For any $\epsilon>0$ there is a $\nu>0$ satisfying the following
condition: For any $n\in\N$, any pair $a,b\in (M_n)_{sa}$ with
$\|b\|\leq1$ and $\|[a,b]\|<\nu$, and any $t\in\R$ there exists a
projection $p\in M_n$ such that
  \BE
  E_a[t+1/4,\infty)&\leq& p\leq E_a(t-1/4,\infty),\\
  \|[a,p]\|&<&\epsilon,\\
  \|[b,p]\|&<&\epsilon,
  \EE
where $E_a$ denotes the spectral measure associated with $a$.
\end{lem}
\begin{pf}
Let $f$ be a $C^{\infty}$-function on $\R$ such that
  $$
  f(t)=\left\{\begin{array}{cl}0&t\leq -1/4\\
                      1& t\geq 1/4 \end{array}\right.
  $$
and $f(t)\approx 2t+1/2,\ 0<f(t)<1$ for $t\in (-1/4,1/4)$. Define
a function $g_N$ on $\R$ for a large $N$ by
  $$
  g_N(t)=\min\{f(t),f(\sqrt{N}-t/\sqrt{N})\}.
  $$
The function $g_N$ is $C^{\infty}$ if $N-\sqrt{N}/4>1/4$ and
satisfies that
  $$
  g_N(t)=\left\{ \begin{array}{cl}1& t\in [1/4,N-\sqrt{N}/4]\\
      0& t\leq -1/4\ {\rm or}\ t\geq N+\sqrt{N}/4\end{array}\right..
  $$
If $N-\sqrt{N}/4\geq \|a\|$, we have that
  $$
  f(a)=g_N(a)=\int\hat{g_N}(t)e^{ita}dt,
  $$
where
  $$
  \hat{g_N}(t)=\frac{1}{2\pi}\int g_N(s)e^{-its}ds.
  $$
Since
  $$
  \|[b,f(a)]\|
  =\int \hat{g_N}(t)[b,e^{ita}]dt=\int\hat{g_N}(t)\int_0^t
               e^{i(t-s)a}[b,ia]e^{isa}dsdt,
  $$
we have that
  $$
  \|[b,f(a)\|\leq \int |\hat{g_N}(t)t|dt\cdot \|[b,a]\|.
  $$
Since
  $$
  it\hat{g_N}(t)=-\frac{1}{2\pi}\int g_N(s)\frac{d}{ds}e^{-its}ds
  =\frac{1}{2\pi}\int g_N'(s)e^{-its},
  $$
it follows for $t\neq0$ that:
  \BE
  \lim_{N\rightarrow\infty}it\hat{g_N}(t)
  &=&\frac{1}{2\pi}\int f'(s)e^{-its}ds
   -\lim \frac{1}{2\pi\sqrt{N}}\int f'(\sqrt{N}-s/\sqrt{N})
        e^{-its}ds\\
    &=& \frac{1}{2\pi}\int f'(s)e^{-its}ds\\
    &=&\hat{f'}(t).
   \EE
Since the above convergence can be estimated by
  $$
  \frac{1}{2\pi\sqrt{N}}\int
  f'(\sqrt{N}-s/\sqrt{N})e^{-its}ds=\frac{e^{-iNt}}{2\pi}\int
  f'(u)e^{i\sqrt{N}tu}du=e^{-iNt}\hat{f'}(-\sqrt{N}t),
  $$
we obtain that
  $$
  \|[b,f(a)]\|\leq C\|[b,a]\|,
  $$
where
   $$
   C=\int |\hat{f'}(t)|dt.
   $$
If $\|[a,b]\|$ is small enough, then $\|[b,f(a)]\|$ is so small
with $\|f(a)\|\leq1$ and $\|b\|\leq1$ that H. Lin's result is
applicable to the pair $b,c=f(a)$. Thus we obtain $b_1,c_1\in
(M_n)_{sa}$ such that
  $$
  \|b-b_1\|\approx0,\ \|c-c_1\|\approx0,\ [b_1,c_1]=0.
  $$
Let $q$ be the spectral projection of $c_1$ corresponding to
$(1/2,\infty)$. Since $\|c-c_1\|\approx0$, and the spectral
projection of $c$ corresponding to $(0,\infty)$ (resp.
$[1,\infty)$) is $E_a(-1/4,\infty)$ (resp. $E_a[1/4,\infty)$), we
have that
  \BE
  E_a(-1/4,\infty)q&\approx& q,\\
  E_a[1/4,\infty)q&\approx& E_a[1/4,\infty),
  \EE
where the approximation depends only on $\|c-c_1\|$, which in turn
depends only on $\|[a,b]\|$. Hence in particular $E_a(-1/4,1/4)$
almost commutes with $q$. By functional calculus we construct a
projection $q_0$ from $E_a(-1/4,1/4)qE_a(-1/4,1/4)$ and set
$p=q_0+E_a[1/4,\infty)$, which is close to $q$, dominates
$E_a[1/4,\infty)$ and is dominated by $E_a(-1/4,\infty)$. Since
$[p,a]=[p,E_a(-1/4,1/4)a]
=[p-q,E_a(-1/4,1/4)a]+[q,E_a(-1/4,1/4)(a-f(a)/2+1/4)]$, we obtain
that $\|[p,a]\|\leq
2\|p-q\|+2\sup_{t\in(-1/4,1/4)}|t-f(t)/2+1/4|$. Since
$[p,b]=[p-q,b]+[q,b]=[p-q,b]+[q,b-b_1]+[q,b_1]$, we obtain that
$\|[p,b]\|\leq 2\|b\|\|p-q\|+2\|b-b_1\|$. Hence $p$ is the desired
projection for $t=0$.  We can apply this argument to the pair
$a-t1,b$ to obtain the desired projection $p$ for $t\in\R$.

\end{pf}

\begin{lem}\label{C3}
For any $\epsilon>0$ there exists a $\nu>0$ satisfying the
following condition:  For any $n \in\N$, any pair $a,b\in
(M_n)_{sa}$ with $\|b\|\leq1$ and $\|[a,b]\|<\nu$ there is a
family $\{p_k:\ k\in \Z\}$ of projections in $M_n$ such that
  \BE
   &&[E_a(j-1/4,j+1/4),p_k]=0,\ \ j,k\in\N,\\
   &&E_a[k+1/4,k+3/4]\leq p_k\leq E_a(k-1/4,k+5/4),\\
   &&\|[a,p_k]\|<\epsilon,\\
   &&\|[b,p_k]\|<\epsilon,\\
   &&\sum_{k}p_k=1,
   \EE
where $p_k=0$ except for a finite number of $k$.
\end{lem}
\begin{pf}
By the previous lemma we choose a $\nu>0$ such that for a pair
$a,b$ as above, there are projections $e_k,  \ k\in\Z$ such that
  \BE
  &&E_a[k+1/4,\infty)\leq e_k\leq E_a(k-1/4,\infty),\\
  &&\|[a,e_k]\|<\epsilon/2,\\
  &&\|[b,e_k]\|<\epsilon/2.
  \EE
Then we set
  $$
  p_k=e_k(1-e_{k+1})=e_k-e_{k+1}.
  $$
Then $\{p_k\}$ is a family of projections with $\sum_kp_k=1$.
Since
  $$E_a(-\infty,k+3/4]\leq 1-e_{k+1}\leq E_a(-\infty,k+5/4),
  $$
we see that $\{p_k\}$ satisfies the required conditions.
\end{pf}

\noindent {\em Proof of Theorem \ref{C}}

By Lemma \ref{C1} we may assume that we are given a pair $a,b\in
(M_n)_{sa}$ such that $\|b\|\leq 1,\ \|[a,b]\|<\nu$, and
$E_a(-\infty,t-1/4]bE_a[t+1/4,\infty)=0$ for any $t\in\R$, where
$\nu>0$ is given in the previous lemma. Choosing the projections
$\{p_k\}$ given there, we claim that
  \BE
  &&\|a-\sum_kp_kap_k\|<4\epsilon,\\
  &&\|b-\sum_kp_kbp_k\|<4\epsilon.
  \EE
To prove this note that if $|i-j|>1$ then $p_iap_j=0=p_ibp_j$.
Since
  $$
  a-\sum_kp_kap_k=\sum_kp_kap_{k+1}+\sum_kp_{k+1}ap_k
  =\sum_k[p_k,a]p_{k+1}+\sum_kp_{k+1}[a,p_k],
  $$
and
  $$\|\sum_{k}[p_{2k},a]p_{2k+1}\|^2=\|\sum_k[p_{2k},a]p_{2k+1}
  [a,p_{2k}]\|
  =\sup_k\|[p_{2k},a]p_{2k+1}[a,p_{2k}]\| <\epsilon^2,
  $$
and similar computations hold for three other sums and for $b$, we
get the above assertions. We then apply H. Lin's result \cite{L}
to each pair $p_kap_k,p_kbp_k$ which satisfies
  $$
  \|[p_kap_k,p_kbp_k]\|
  \leq \|p_k[a,p_k]bp_k\|+\|p_k[a,b]p_k\|+\|p_k[b,p_k]ap_k\|
   <2\epsilon+\nu.
  $$
Assuming that $2\epsilon+\nu$ is sufficiently small, we obtain a
pair $a_k,b_k$ in $(p_kM_np_k)_{sa}$ such that
  $$
  p_kap_k\approx a_k,\ p_kbp_k\approx b_k,\ [a_k,b_k]=0.
  $$
We set $a'=\sum_ka_k$ and $b'=\sum_kb_k$. Then it follows that
$[a',b']=0$ and $a\approx a',\ b\approx b'$ because of the
inequality
  $$
  \|a-a'\|\leq \|a-\sum_kp_kap_k\|+\sup_k \|p_kap_k-a_k\|
  $$
and a similar inequality for $b,b'$. This completes the proof.
\medskip

For a flow $\alpha$ of a unital simple AF algebra we denote by
$\delta_{\alpha}$ the generator of $\alpha$ as before. We
introduce the following condition on $\alpha$, which we may
express by saying that {\em the almost fixed point algebra for
$\alpha$ has real rank zero.}

\medskip
\noindent
{\bf Condition F0:}  For any $\epsilon>0$ there exists a
$\nu>0$ satisfying the following condition: If $h=h^*\in
\D(\delta_{\alpha})$ satisfies that $\|h\|\leq 1$ and
$\|\delta_{\alpha}(h)\|<\nu$ there exists a pair $k=k^*\in
\D(\delta_{\alpha})$ and $b=b^*\in A$ such that
$\|h-k\|<\epsilon,\ \|b\|<\epsilon, \
(\delta_{\alpha}+\ad\,ib)(k)=0$, and $\Sp(k)$ is finite.

\medskip
In the above condition let $C$ be the (finite-dimensional)
*-subalgebra generated by $k$. Then $h$ is approximated by an
element of $C$ within distance $\epsilon$ and
$\|\delta_{\alpha}|C\|<2\epsilon$.

We recall  from \cite[Proposition~3.1]{Kis99} that a flow
$\alpha$ is a cocycle perturbation of an AF flow if and only if the
domain
$\D(\delta_{\alpha})$ contains a canonical AF masa. (A maximal
abelian AF
$C^*$-subalgebra  $C$ of a AF $C^*$-algebra $A$ is called {\em
canonical} if there is an increasing sequence $(A_n)$ of
finite-dimensional *-subalgebras of
$A$ with dense union such that $C\cap A_{n}\cap A_{n-1}'$ is
maximal abelian in $A_n\cap A_{n-1}'$ for each $n$ with $A_0=0$.)

\begin{theo}\label{C4}
Let $\alpha$ be a flow of a non type I simple AF $C^*$-algebra. If
$\D(\delta_{\alpha})$ contains a canonical AF masa, then the above
condition F0 is satisfied, i.e., the almost fixed point algebra
has real rank zero.
\end{theo}
\begin{pf}
Let $\epsilon>0$. We choose a $\nu>0$ as in Theorem \ref{C}.

Let $h=h^*\in \D(\delta_{\alpha})$ be such that $\|h\|\leq1$ and
$\|\delta_{\alpha}(h)\|<\nu$. There exists a $c=c^*\in A$ such
that $\|c\|<\min\{(\nu-\|\delta_{\alpha}(h)\|)/2,\epsilon\}$ and
$\delta_{\alpha}+\ad\,ic$ generates an AF flow. Explicitly let
$\{A_n\}$ be an increasing sequence of finite-dimensional
*subalgebras of $A$ with dense union such that $A_n\subset
\D(\delta_{\alpha})$ and $(\delta_{\alpha}+\ad\,ic)(A_n)\subset
A_n$ for each $n$. There exists a sequence $\{h_n\}$ such that
$h_n=h_n^*\in A_n,\ \|h_n\|\leq1,\ \|h_n-h\|\ra0$, and
$\|\delta_{\alpha}(h-h_n)\|\ra0$. Since
$\|(\delta_{\alpha}+\ad\,ic)(h)\|<\nu$, we have an $n$,
$h_0=h_0^*\in A_n$, and $a=a^*\in A_n$ such that $\|h_0\|\leq1$,
$\|h-h_0\|<\epsilon$, $\|(\delta_{\alpha}+\ad\,ic)(h_0)\|<\nu$,
and $(\delta_{\alpha}+\ad\,ic)|A_n=\ad\,ia|A_n$. Since $A_n$ is a
finite direct sum of matrix algebras, Theorem \ref{C} is
applicable to the pair $a,h_0$. Thus there exists a pair
$a_1,h_1\in (A_n)_{sa}$ such that $\|a-a_1\|<\epsilon,\
\|h_0-h_1\|<\epsilon$, and $[a_1,h_1]=0$. Let $b=a_1-a+c$. Then we
have that $\|h-h_1\|<2\epsilon,\ \|b\|<2\epsilon,\
(\delta_{\alpha}+\ad\,ib)(h_1)=0$, and $\Sp(h_1)$ is finite.
\end{pf}

In the special case that $\alpha$ is periodic, the fact that the
almost fixed point algebra has real rank zero simply means that the
fixed point algebra has real rank zero:

\begin{prop}\label{C5}
Let $A$ a non type I simple AF $C^*$-algebra and $\alpha$ a
periodic flow of $A$. Then the following conditions are
equivalent:
\begin{enumerate}
\item Condition F0 holds.
\item The fixed point algebra
$A^{\alpha}=\{a\in A\ | \ \alpha_t(a)=a\}$ has real rank zero.
\end{enumerate}
\end{prop}
\begin{pf}
We may suppose that $\alpha_1=\id$. Suppose (1); we have to show
that $\{h\in A_{sa}^{\alpha}\ |\ \Sp(h)\ {\rm is\ finite}\}$ is
dense in $A_{sa}^{\alpha}$ \cite{BP}. Let $h=h^*\in A^{\alpha}$,
$\epsilon>0$, and $n\in\N$. There exist an $h_1\in
{\D}(\delta_{\alpha})_{sa}$ and  $b\in A_{sa}$ such that
$\|h-h_1\|<\epsilon,\ \|b\|<\epsilon,\
(\delta_{\alpha}+\ad\,ib)(h_1)=0$, and $\Sp(h_1)$ is finite. We
approximate $h_1$ by an element $h_2=\sum_{k=-n}^n (k/n)p_k$ in
the *-subalgebra generated by $h_1$, where $(p_k)$ is a mutually
orthogonal family of projections. We may assume that
$\|h_1-h_2\|\leq 1/n$ and hence that $\|h-h_2\|<\epsilon+1/n$.
Note that we still have that $(\delta_{\alpha}+\ad\,ib)(h_2)=0$.
Since $\|\alpha_t(p_k)-p_k\|\leq |t|\|\delta(p_k)\|<2|t|\epsilon$,
we have that
  $$
  \|\int_0^1\alpha_t(p_k)-p_k\|<\epsilon
  $$
for $k=-n,-n+1,\ldots,n$. If $\epsilon$ is sufficiently small,
then by functional calculus we inductively define a projection
$q_k\in A^{\alpha}$ from
$(1-\sum_{j=-n}^{k-1}q_j)\int\alpha_tp_kdt(1-\sum_{j=-n}^{k-1}q_j)$,
which belongs to $A^{\alpha}$, such that $q_k\approx p_k$ and
$q_k$ is orthogonal to $\sum_{j=-n}^{k-1}q_j$. Then
$h_3=\sum_{k=-n}^n(k/n)q_k\approx \sum_{k=-n}^n(k/n)p_k=h_2$,
where the approximation is of the order of $\epsilon$ times some
function of $n$. Since $h_3\in A^{\alpha}$, we reach the
conclusion by choosing $\epsilon>0$ sufficiently small.

The converse implication is easy to show. \end{pf}

If $\alpha$ is not periodic, we can still re-formulate Condition
F0 as follows, further justifying the terminology that the almost
fixed point algebra has real rank zero. We denote by
$\ell^{\infty}$ the
\cstar\ of bounded sequences in $A$ and by $c_0$ the closed ideal
of
$\ell^{\infty}$ consisting of sequences converging to zero. Then we
set
$A^{\infty}$ to be the quotient $\ell^{\infty}/c_0$. The flow
$\alpha$ on $A$ induces a flow $\ol{\alpha}$ on $\ell^{\infty}$ by
$\ol{\alpha}_t(x)=(\alpha_t(x_n))$ for $x=(x_n)$. But since
$\ol{\alpha}$ is not strongly continuous (if $\alpha$ is not
uniformly continuous), we choose the $C^*$-subalgebra
$\ell_{\alpha}^{\infty}$ consisting of $x\in\ell^{\infty}$ with
$t\mapsto \ol{\alpha}_t(x)$ continuous. Since
$\ell_{\alpha}^{\infty}\supset c_0$ and $c_0$ is
$\ol{\alpha}$-invariant, $\ol{\alpha}$ induces a (strongly
continuous) flow on the quotient
$A_{\alpha}^{\infty}=\ell_{\alpha}^{\infty}/c_0$, which will also
be denoted by $\alpha$. Note that $A_{\alpha}^{\infty}$ is
inseparable even if $A$ is separable. See \cite{Kis96}.

\begin{prop}\label{C51}
Let $A$ be a \cstar\ and $\alpha$ a flow of $A$. Then the
following conditions are equivalent:
  \begin{enumerate}
  \item Condition F0 holds.
  \item The fixed point algebra $(A_{\alpha}^{\infty})^{\alpha}$
  has real rank zero.
  \end{enumerate}
\end{prop}
\begin{pf}
Suppose (1) and let $h\in (A_{\alpha}^{\infty})_{sa}^{\alpha}$. We
take a representative $(h_n)\in\ell_{\alpha}^{\infty}$ of $h$ such
that $h_n^*=h_n$ for all $n$. Taking a non-negative $C^{\infty}$
function $f$ with integral 1, we may replace each $h_n$ by $\int
f(t)\alpha_t(h_n)dt$. Thus we can assume that
$h_n\in\D(\delta_{\alpha})$ and $\|\delta_{\alpha}(h_n)\|\ra0$.
Then for any $\epsilon>0$ there exists a sequence of pairs $k_n\in
\D(\delta_{\alpha})_{sa}$ and $b_n\in A_{sa}$ such that
$\|h_n-k_n\|<\epsilon,\ \|b_n\|\ra0,
(\delta_{\alpha}+\ad\,ib_n)(k_n)=0$, and $\Sp(k_n)$ is finite and
independent of $n$. Hence $k=(k_n)+c_0\in A_{\alpha}^{\infty}$
satisfies that $\|h-k\|\leq \epsilon,\ \delta_{\alpha}(k)=0$, and
$\Sp(k)$ is finite. This shows that
$(A_{\alpha}^{\infty})^{\alpha}$ has real rank zero \cite{BP}.

Suppose (2). If Condition F0 does not hold, we find an
$\epsilon>0$ and a sequence $(h_n)$ in $\D(\delta_{\alpha})_{sa}$
such that $\|h_n\|=1,\ \|\delta_{\alpha}(h_n)\|\ra0$, and such
that if $k\in \D(\delta_{\alpha})_{sa}$ and $b\in A_{sa}$ satisfy
that $\|h-k\|<\epsilon,\ \|b\|<\epsilon$, and $\Sp(k)$ is finite,
then $(\delta_{\alpha}+\ad\,ib)(k)\neq0$. Since $h=(h_n)+c_0\in
A_{\alpha}^{\infty}$ belongs to $(A_{\alpha}^{\infty})^{\alpha}$,
we have a $k\in (A_{\alpha}^{\infty})_{sa}^{\alpha}$ such that
$\|h-k\|<\epsilon$ and $\Sp(k)$ is finite. By choosing an
appropriate representative (consisting of projections) for each
minimal spectral projection of $k$, we find a representative
$(k_n)$ of $k$ such that $k_n^*=k_n,\ \Sp(k_n)=\Sp(k)$, and
$\|\delta_{\alpha}(k_n)\|\ra0$. This is a contradiction. \end{pf}

We recall here a condition on a flow $\alpha$ considered in
\cite{Kis99}.

\medskip
\noindent
{\bf Condition F1:} For any $\epsilon>0$ there exists a
$\nu>0$ satisfying the following condition: If
$u\in\D(\delta_{\alpha})$ is a unitary with
$\|\delta_{\alpha}(u)\|<\nu$ there is a continuous path $(u_t)$ of
unitaries in $A$ such that $u_0=1,\ u_1=u,\ u_t\in
\D(\delta_{\alpha})$, and $\|\delta_{\alpha}(u_t)\|<\epsilon$ for
$t\in[0,1]$.

\medskip
In the above condition we can choose the path $(u_t)$ to be
continuous in the Banach *-algebra $\D(\delta_{\alpha})$. We
express this condition by saying that {\em the almost fixed point
algebra for $\alpha$ has trivial $K_1$}. What we have shown in
\cite{Kis99} is that if $\alpha$ is an inner perturbation of an AF
flow then the above condition holds. Actually by using the full
strength of Lemma 5.1 of \cite{BBEK}, one can show that the
following stronger condition holds:

\medskip
\noindent
{\bf Condition F1':}  For any $\epsilon>0$ there exists
a $\nu>0$  satisfying the following condition: If
$u\in\D(\delta_{\alpha})$ is a unitary with
$\|\delta_{\alpha}(u)\|<\nu$ there is a rectifiable path $(u_t)$
of unitaries in $A$ such that $u_0=1,\ u_1=u,\ u_t\in
\D(\delta_{\alpha})$, $\|\delta_{\alpha}(u_t)\|<\epsilon$ for
$t\in[0,1]$, and the length of $(u_t)$ is bounded by $C$, where
$C$ is a universal constant (smaller than $3\pi+\varepsilon$ for
example).

\medskip
Then one can show the following:
\begin{prop}
Let $A$ be a unital  \cstar\ and $\alpha$ a flow of $A$. Then the
following conditions are equivalent:
  \begin{enumerate}
  \item Condition F1' holds.
  \item The unitary group of the fixed point algebra
  $(A_{\alpha}^{\infty})^{\alpha}$
  is path-wise connected; moreover any unitary is connected to 1
  by a continuous path of unitaries whose length is bounded by
  a universal constant.
  \end{enumerate}
\end{prop}
We will leave the proof to the reader.

\begin{rem}\label{C6}
If $A$ is a unital simple AF $C^*$-algebra, one can construct a
periodic flow $\alpha$ of $A$, by using the general classification
theory of locally representable actions \cite{BEEK}, such that the
almost fixed point algebra for $\alpha$ has real rank zero but
does not have trivial $K_1$.
\end{rem}

\begin{prop}\label{C7}
Let $A$ be a unital simple AF \cstar. Then there exists a flow
$\alpha$ of $A$ such that $\D(\delta_{\alpha})$ is AF and the
almost fixed point algebra for $\alpha$ does not have real rank
zero but has trivial $K_1$ (i.e., F0 holds but not F1).
\end{prop}
\begin{pf}
We shall use a construction used in the proof of 2.1 of
\cite{Kis99}. Let $(A_n)$ be an increasing sequence of
finite-dimensional *-subalgebras of $A$ such that
$A=\ol{\cup_nA_n}$
and let $A_n=\oplus_{j=1}^{k_n}A_{nj}$ be the direct sum
decomposition of $A_n$ into full matrix algebras $A_{nj}$. Since
$K_0(A_n)\cong\Z^{k_n}$, we obtain a sequence of $K_0$ groups:
  $$
  \Z^{k_1}\stackrel{\chi_1}{\ra}\Z^{k_2}\stackrel{\chi_2}{\ra}
     \cdots,
  $$
where $\chi_n$ is the positive map of $K_0(A_n)=\Z^{k_n}$ into
$K_0(A_{n+1})=\Z^{k_{n+1}}$ induced by the embedding $A_n\subset
A_{n+1}$. Since $K_0(A)$ is a simple dimension group different from
$\Z$, we may assume that $\min_{ij}\chi_n(i,j)\ra \infty$ as
$n\ra\infty$.

By using $(A_n)$ we will express $A$ as an inductive limit of
\cstars\ $A_n\otimes C[0,1]$. First we define a homomorphism
$\varphi_{n,ij}$ of $A_{nj}\otimes C[0,1]$ into $A_{nj}\otimes
M_{\chi_n(i,j)}\otimes C[0,1]$ as follows: If $i=j=1$ then
  $$
  \varphi_{n,11}(x)(t)=x(t)\oplus
  \oplus_{\ell=0}^{\chi_n(1,1)-2}x(\frac{t+\ell}{\chi_n(1,1)-1}),
  $$
otherwise
  $$
  \varphi_{n,ij}(x)(t)
  = \oplus_{\ell=0}^{\chi_n(1,1)-1}x(\frac{t+\ell}{\chi_n(1,1)}).
  $$
Especially $\varphi_{n,ij}(x)$ is of diagonal form in the matrix
algebra over $A_{nj}\otimes C[0,1]$. Then embedding
  $$
  \oplus_{j=1}^{k_n}A_{nj}\otimes M_{\chi_n(i,j)}\otimes C[0,1]
  $$
into $A_{n+1,i}\otimes C[0,1]$, $(\varphi_{n,ij})$ defines an
injective homomorphism $\varphi_n:\,A_n\otimes C[0,1]\ra
A_{n+1}\otimes C[0,1]$. Then it follows that the inductive limit
\cstar\ of $(A_n\otimes C[0,1],\varphi_n)$ is isomorphic to the
original $A$; we have thus expressed $A$ as $\ol{\cup_nB_n}$ where
$B_n=A_n\otimes C[0,1]\subset B_{n+1}$ \cite{Ell1}.

We will define a flow or one-parameter automorphism group $\alpha$
of $A$ such that $\alpha_t(B_n)=B_n$ and $\alpha_t\big|_{B_n}$ is
inner, i.e.,
$\alpha$ is locally representable for the sequence $(B_n)$. First
we define a sequence $(H_n)$ with self-adjoint $H_n\in A_n\otimes
1\subset B_n$ inductively. Let $H_1\in A_1\otimes 1\subset B_1$ and
let
$H_{n}=H_{n-1}+\sum_i\sum_j h_{n,ij}$, where
  $$
  h_{n,ij}^*=h_{n,ij}\in 1\otimes M_{\chi_{n-1}(i,j)}\otimes 1
  \subset A_{n-1,j}\otimes M_{\chi_{n-1}(i,j)}\otimes 1\subset B_n.
  $$
We define $\alpha_t|B_n$ by $\Ad\,e^{itH_n}|B_n$. Since
$\alpha_t|B_n=\Ad\,e^{itH_{n+1}}|B_n$ from the definition of
$H_{n+1}$, $(\alpha_t|B_n)$ defines a flow $\alpha$ of $A$.

We fix $H_1$ and $h_{nij}$ in the following way:
$\|h_{nij}\|\leq1/2$ except for $h_{n11}$ which is defined by
  $$
  h_{n11}
  =1\oplus0\oplus\cdots\oplus\in1\otimes M_{\chi_{n-1}(1,1)}\otimes1
  \subset A_{n1}\otimes C[0,1].
  $$
We will show that the $\alpha$ defined this way has the desired
properties.

Let $x$ be the identity function on the interval $[0,1]$ and let
$x_n=1\otimes x\in 1\otimes C[0,1]\subset B_n$. To show that
$\D(\delta_{\alpha})$ is AF, it suffices to show that for each
$x_n$, there exists a sequence $(h_m)_{m>n}$ such that
$h_m=h_m^*\in B_m$, $\Sp(h_m)$ is finite, and
$\|x_n-h_m\|_{\delta_{\alpha}}\ra0$ as $m\ra\infty$. For a
sufficiently large $m>n$, the image $\varphi_{mn}(x_n)$ of $x_n$
in $B_m=A_m\otimes C[0,1]$ is almost constant as a function (into
the diagonal matrices in $A_m\cap A_n'$) on $[0,1]$ except for one
component, which is $x$ and appears through the first component of
$\varphi_{k11}$ for $n\leq k<m$. We will approximate this
component $x$ by a self-adjoint element with finite spectrum by
using the part appearing through the components of $\varphi_{k11}$
other than the first; they are the direct sum of
$M=\Pi_{k=n}^{m-1}(\chi_k(1,1)-1)$ components
$x(\frac{t+\ell}{M}),\ \ell=0,1,\ldots,M-1$. There is a standard
procedure to approximate the sum of these $M+1$ components by a
self-adjoint element $k$ with finite spectrum \cite{BBEK}. Since
$H_m-H_n$ is $m-n$ on the support projection of $x$ and 0 on the
support projections of the other components, the
$\|\,\cdot\,\|_{\delta_{\alpha}}$ norm of
$k$ is of the order of $\frac{m-n}{M}\approx0$. (All the spectral
projections of $k$ are just constant at each point of $[0,1]$
perhaps except for a pair of projections, whose eigen-values are
different only by the order of $1/M$, and which are of the form:
  $$
  \left(\begin{array}{cc}\cos^2\theta& \cos\theta\sin\theta\\
     \cos\theta\sin\theta& \sin^2\theta\end{array}\right),\ \
   \left(\begin{array}{cc}\sin^2\theta& -\cos\theta\sin\theta\\
     -\cos\theta\sin\theta& \cos^2\theta\end{array}\right)
  $$
in the space spanned by the support projection of $x$ and one of
the support projections of the other $M$ components, where
$\theta$ is a function in $t\in[0,1]$ which changes from $-\pi/2$
to $\pi/2$ quickly near the point in problem. This implies that
$\|\delta_{\alpha}(k)\|\approx \frac{m-n}{M}$ and
$\|x_n-k\|\approx 1/M$ for the parts of  $k,\,x_n-k$ in question.)
This concludes the proof that $\D(\delta_{\alpha})$ is AF.

Suppose that for any $\epsilon>0$ there exists a pair of
self-adjoint elements $h,b\in A$ such that $\|h\|\leq1,\
\|b\|<\epsilon,\ \|x_1-h\|<\epsilon$,
$(\delta_{\alpha}+\ad\,ib)(h)=0$, and $\Sp(h)$ is finite, where
$x_1$ is the element of $B_1$ defined above.  Since $\cup_mB_m$ is
dense in $\D(\delta_{\alpha})$, we may suppose that $h\in B_m$ for
some $m$. The image $\varphi_{m1}(x_1)$ in $B_m\cap A_1'$ is
diagonal and there is a component $x$, whose (one-dimensional)
support projection will be denoted by $Q$. Let
$h=\sum_i\lambda_ip_i$ be the spectral decomposition of $h$ and
define a function $\theta_i$ by $\theta_i(t)=Qp_i(t)Q$. Then we
have that
  $$
  |t-\sum_i\lambda_i\theta_i(t)|<\epsilon,\ \ \ t\in[0,1].
  $$
Since
  $$
  \frac{1}{2}\sum_{\lambda_i>1/2}\theta_i(0)
  <\sum\lambda_i\theta_i(0)<\epsilon,
  $$
we obtain that
  $$
  \sum_{\lambda_i>1/2}\theta_i(0)<2\epsilon.
  $$
Since
  $$
  1-\epsilon<\sum\lambda_i\theta_i(1)
  <\frac{1}{2}\sum_{\lambda\leq1/2}\theta_i(1)
  +\sum_{\lambda_i>1/2}\theta_i(1)
  =\frac{1}{2}+\frac{1}{2}\sum_{\lambda_i>1/2}\theta_i(1),
  $$
we get
  $$
  \sum_{\lambda_i>1/2}\theta_i(1)>1-2\epsilon.
  $$
Thus the projection $p$ defined by
  $$
  p=\sum_{\lambda_i>1/2}p_i
  $$
satisfies that $\|Qp(0)Q\|<2\epsilon$ and
$\|Qp(1)Q\|>1-2\epsilon$. If $\epsilon<1/4$, there must be a point
$t\in [0,1]$ such that $\|Qp(t)Q\|=1/2$. Then since
$Qp(t)(1-Q)p(t)Q+Qp(t)Qp(t)Q=Qp(t)Q$, we have that
$\|Qp(t)(1-Q)\|=1/2$. Since $(H_m-H_1)Q=(m-1)Q$ and
$\|(H_m-H_1)(1-Q)\|\leq m-3/2$, we get that
$\|\delta_{\alpha}(Qp(1-Q))\|=\|Q\delta_{\alpha}(p)(1-Q)\|\geq1/4$.
But since $(\delta_{\alpha}+\ad\,ib)(h)=0$, we had that
$\|\delta_{\alpha}(p)\|\leq 2\|b\|<2\epsilon$. For a small
$\epsilon>0$ this is a contradiction. Thus we obtain that the
almost fixed point algebra does not have real rank zero.

Let $u$ be a unitary in $\D(\delta_{\alpha})$ such that
$\delta_{\alpha}(u)\approx0$. Since $\cup_mB_m$ is dense in
$\D(\delta_{\alpha})$, we may suppose that $u\in B_m=A_m\otimes
C[0,1]$. Since $H_m\in A_m\otimes 1$, the condition
$\delta_{\alpha}(u)\approx0$ implies that $\|[u(t),H_m]\|\approx0$
for all $t\in [0,1]$. Define a continuous path $(u_s)$ of
unitaries in $B_m$ by $u_s(t)=u((1-s)t)$. This path runs from $u$
to the constant function $u_1:t\mapsto u(0)$ with the estimate
$\|\delta_{\alpha}(u_s)\|\leq \|\delta_{\alpha}(u)\|$. By 4.1 of
\cite{Kis99}, there is a continuous path $(v_s)$ of
unitaies in $A_m$ from $u(0)$ to 1 such that $[v_s,H_m]\approx0$.
 This concludes the proof that the
almost fixed point algebra has trivial $K_1$.
\end{pf}

\section{The CAR algebra}
\setcounter{theo}{0} \ncm{\A}{{\cal A}} \ncm{\Hil}{{\cal H}}
\ncm{\M}{{\cal M}} Let $A=\A(\Hil)$ be the CAR algebra over an
infinite-dimensional separable Hilbert space $\Hil$; we denote by
$a^*$ the canonical linear isometric map of $\Hil$ into the
creation operators in $A$, \cite[Section~5.2.2.1]{BR2}. Note that
$A$, as a
\cstar, is isomorphic to the UHF algebra of type
$2^{\infty}$. When $U$ is a one-parameter unitary group on $\Hil$,
we define a flow $\alpha$ of $A$ by
  $$
  \alpha_t(a^*(\xi))=a^*(U_t\xi),\ \ \ \xi\in\Hil,
  $$
which will be called the quasi-free flow induced by $U$. If we
denote by $H$ the generator of $U$, i.e., $U_t=e^{itH}$, the
generator $\delta_{\alpha}$ of $\alpha$ satisfies that
  $$
  \delta_{\alpha}(a^*(\xi))=ia^*(H\xi),\ \ \ \xi\in\D(H)
  $$
and the *-subalgebra generated by $a^*(\xi),\ \xi\in \D(H)$ is
dense in the Banach *-algebra $\D(\delta_{\alpha})$. If $H$ is
diagonal, i.e., has a complete orthonormal family of eigenvectors,
then $\alpha$ is an AF flow; moreover it is of {\em of pure
product type} in the sense that $(A,\alpha)$ is isomorphic to
$(M_{2^{\infty}},\beta)$, where $\beta$ is given as
  $$
\beta_t=\otimes_{n=1}^{\infty}\Ad\,\left(
\begin{array}{cc}e^{i\lambda_nt}&0\\
 0&1 \end{array}\right),
  $$
where $\{\lambda_n, n\in \Z\}$ are the eigenvalues of $H$. If $H$
is not diagonal,
$\alpha$ acts on a part of
$A$ {\em in an asymptotically abelian way}; so we can conclude that
$\alpha$ is not an AF flow.  See \cite[8, 18]{BR2} for
details.

\begin{prop}\label{D}
If $\alpha$ is a quasi-free flow of the CAR algebra $A=\A(\Hil)$,
then the almost fixed point algebra for $\alpha$ has trivial
$K_1$. \end{prop}
\begin{pf}
We use the notation given before this proposition and let $E$ be
the spectral measure of $H$. Let $\epsilon>0$ and let
$u\in\D(\delta_{\alpha})$ be a unitary such that
$\|\delta_{\alpha}(u)\|<\epsilon$. Since the *-subalgebra ${\cal
P}$ generated by
  $$
  a^*(\xi),\ \ \ \xi\in\bigcup E[-n,n]\Hil
  $$
is dense in $\D(\delta_{\alpha})$, we can approximate $u$ by
$x\in{\cal P}$. Let $\M$ be the (abelian) von Neumann algebra
generated by $U_t=e^{itH},\ t\in\R$. We may approximate $u$ by $x$
in a *-subalgebra ${\cal P}_1$ generated by
$a^*(\xi_1),a^*(\xi_2),\ldots,a^*(\xi_n)$, where all $\xi_i\in
E[-N,N]\Hil$ for some $N$. We may further impose the following
conditions on $\xi_1,\ldots,\xi_n$:
  \begin{enumerate}
  \item $\|\xi_i\|=1$ for all $i$.
  \item For $i\neq j$, $\ol{\M{\xi_i}}\,\bot\,\ol{\M{\xi_j}}$.
  \item Denote by $S_i$ the smallest closed subset of $\R$ such
that
  $E(S_i)\xi=\xi$. Then either $S_i$ is a singleton or an infinite
set.
  \end{enumerate}
The condition 1 is trivial and the condition 3 is easy to obtain.
To make sure the condition 2 holds we may argue as follows. Starting
with $\xi_1,\ldots,\xi_n$ let $e_1$ be the projection onto
$\ol{\M\xi_1}$. Then
$\xi_1'=\xi_1=e_1\xi_1,\xi_2'=e_1x_2,\ldots,\xi_n'=e_1\xi_n$ all
belong to $e_1\Hil$ on which $\M e_1$ is a maximal abelian von
Neumann algebra. Thus there are a finite number of unit vectors
$\eta_{11},\ldots,\eta_{1m}$ in $e_1\Hil$ such that the linear
span of $\eta_{1i}$'s approximately contains all $\xi_j'$ and
$\ol{\M\eta_{1i}}\,\bot\, \ol{\M\eta_{1j}}$ for $i\neq j$.  We
apply the same argument to the remaining (at most $n-1$) elements
$(1-e_1)\xi_2,(1-e_1)\xi_3,\ldots,(1-e_1)\xi_n$ in $(1-e_1)\Hil$
which is left invariant under $\M$. Next, let $e_2$ be the
projection onto $\ol{\M(1-e_1)\xi_2}$ (assuming this is non-zero).
Note that
$e_2\leq 1-e_1$. We find a finite number of unit vectors
$\eta_{2j}$ in $e_2\Hil$ whose linear span approximately contains
$e_2(1-e_1)\xi_2=(1-e_1)\xi_2,e_2(1-e_1)\xi_3
=e_2\xi_3,\ldots,e_2(1-e_1)\xi_n=e_2\xi_n$
such that $\ol{\M\eta_{2i}}\,\bot\,\ol{\M\eta_{2j}}$ for $i\neq
j$. Note that $\ol{\M\eta_{1i}}\,\bot\,\ol{\M\eta_{2j}}$ for all
$i,j$. Repeating this procedure we obtain a finite number of unit
vectors $(\eta_{ij})$ satisfying the condition 2 whose linear span
approximately contains the vectors $\xi_1,\ldots,\xi_n$.

Since the *-algebra ${\cal P}_1$ is isomorphic to $M_{2^n}$ by
\cite[Theorem~5.2.5]{BR2}, we may further assume that $x$ is a
unitary. We express
$x$ as
  $$
  x=\sum_{\mu\nu}a_{\mu\nu}a^*(\mu)a(\nu),
  $$
where $\mu=(\mu_1,\ldots,\mu_{\ell})$ and $\nu$ ranges over the
subsequences of
$(1,2,\ldots,n)$ and $a^*(\mu)$ denotes
  $$a^*(\xi_{\mu_1})a^*(\xi_{\mu_2})\cdots a^*(\xi_{\mu_{\ell}})
  $$
with $a(\nu)=a^*(\nu)^*$. (If $\mu$ is the empty sequence, then
$a^*(\mu)=1$.) Note that the coefficients $a_{\mu\nu}$ are unique;
hence the condition that $x$ is a unitary can be read from
$(a_{\mu\nu})$ only, i.e., if we replace $\xi_1,\ldots,\xi_n$ by a
different orthonormal family $\eta_1,\ldots,\eta_n$ and define $x$
by the same formula with $a^*(\mu)=a^*(\eta_{\mu_1})\cdots
a^*(\eta_{\mu_n})$, then $x$ is still a unitary.

Let
  $$
  \eta_i=H\xi_i-(H\xi_i|\xi_i)\xi_i\;.
  $$
If $\eta_i\neq0$ let $\xi_{i+1/2}=\eta_i/\|\eta_i\|$ and
otherwise let $\xi_{i+1/2}=0$. Let $I$ be the subsequence of
$\{1,3/2,2,\ldots,n+1/2\}$ with $\xi_c\neq0$. Then the vectors
$\xi_c,\ c\in I$ form an orthonormal family. Since
$H\xi_i=(H\xi_i|\xi_i)\xi_i+\|\eta_i\|\xi_{i+1/2}$,
$\delta_{\alpha}(x)$ is of the form
  $$
  \delta_{\alpha}(x)=\sum_{\sigma\tau}b_{\sigma\tau}a^*
(\sigma)a(\tau),
  $$
where $\sigma,\tau$ are subsequences of $I$. Again the norm
$\|\delta_{\alpha}(x)\|$ can be read from $(b_{\sigma\tau})$ only.
Note that $(b_{\sigma\tau})$ depends only on $(a_{\mu\nu})$,
$(H\xi_i|\xi_i)$, and $\|H\xi_i-(H\xi_i|\xi_i)\xi_i\|$.

By using Lemma~4.2 below, if $S_i$ is not a singleton, we will
find a continuous path of $(\xi_{it})_{0\leq t<1}$ of unit vectors
in $\ol{\M\xi_i}\subset \Hil$ such that $\xi_{i0}=\xi_i$,
  \BE
  && (H\xi_{it}|\xi_{it})=(H\xi_i|\xi_i),\\
  &&
\|H\xi_{it}-(H\xi_{it}|\xi_{it})\xi_{it})\|
=\|H\xi_i-(H\xi_i|\xi_i)\xi_i\|,
  \EE
and $\supp\, \xi_{it}$ shrinks to a three point set as $t\ra1$,
where $\supp\, \xi$ is the smallest closed subset $S$ of $\R$ with
$E(S)\xi=\xi$. And we set
$\xi_{i+1/2,t}=c_i(H\xi_{it}-(H\xi_{it}|\xi_{it})\xi_{it})$, where
$c_i$ is a positive normalizing constant. If $S_i$ is a singleton,
we set $\xi_{it}=\xi_i$. By using $\xi_{it},\ 1\leq i\leq n$
instead of $\xi_{i}$, we define $x_t\in A$ by the same formula as
$x$. Then we have that $(x_t)_{0\leq t<1}$ is a continuous family
of unitaries with $x_0=x$ satisfying
  $$
  \delta_{\alpha}(x_t)=\sum_{\sigma\tau}
b_{\sigma\tau}a^*(\sigma)a(\tau),
  $$
where $\sigma, \tau$ are subsequences of $I$ and
$a^*(\sigma),a(\tau)$ are now defined by using $\xi_{ct},\ c\in
I$ instead of
$\xi_c$.  Hence it follows that
$\|\delta_{\alpha}(x_t)\|=\|\delta_{\alpha}(x)\|<\epsilon$.

We will show that for a $t_0$ close to 1 there is a $b=b^*\in A$
such that $\|b\|<\epsilon/2$ and $\delta_{\alpha}+\ad\,ib$ leaves
a finite-dimensional *-subalgebra containing $x_{t_0}$ invariant,
and such that $\|(\delta_{\alpha}+\ad\,ib)x_{t_0}\|$ is sufficiently
small. Then by \cite{BEEK98} we can deform $x_{t_0}$ to 1 in that
*-subalgebra keeping the norm estimate along the path.

Suppose that $t_0$ is sufficiently close to 1. If $S_i$ is a
singleton, we set $\eta_{i1}=\xi_i$; otherwise we choose three
unit vectors $\eta_{i1},\eta_{i2},\eta_{i3}$ in $\ol{\M\xi_i}$
such that $\xi_{it_0}$ is a linear combination of $\eta_{ij}$'s
and $\supp\,\eta_{ij}$ is contained in a sufficiently small
neighborhood of some $s_{ij}\in S_i$, where $s_{i1},s_{i2},s_{i3}$
are all distinct. Let $P_{ij}$ be the projection onto the space
spanned by $\eta_{ij},H\eta_{ij}$ and define an operator $T_{ij}$
such that $T_{ij}=P_{ij}T_{ij}P_{ij}=T_{ij}^*$ and
$T_{ij}\eta_{ij}=(s_{ij}1-H)\eta_{ij}$. Then it follows that the
projections $P_{ij}$ are mutually orthogonal and $\|T_{ij}\|\leq
  2\|(s_{ij}1-H)\eta_{ij}\|$, which is assumed to be very small.
Let
$T_i=\sum_jT_{ij}$ and $T=\sum_i T_i$, where we set $T_i=0$ if
$S_i$ is a singleton.  Then $\|T\|=\sup\,\|T_i\|,\ {\rm
rank}(T)\leq 6n$, and $(H+T)\eta_{ij}=s_{ij}\eta_{ij}$. We may
suppose that ${\rm Tr}(|T|)<\epsilon/2$. Note that the derivation
of $A$ corresponding to $T$ is inner and given as $\ad\,ib$, where
$b=\sum\lambda_ia^*(\zeta_i)a(\zeta_i)$, if $(\zeta_i)$ is a
complete orthonormal set of eigenvectors of $T$ with $(\lambda_i)$
the corresponding eigenvalues; $T\zeta_i=\lambda_i\zeta_i$. If
${\cal P}_2$ denotes the *-algebra generated by $a^*(\eta_{ij})$,
then ${\cal P}_2$ is left invariant under the derivation
corresponding to $H+T$, which is $\delta_{\alpha}+\ad\,ib$. Hence
there is an $h=h^*\in{\cal P}_2$ such that
$(\delta_{\alpha}+\ad\,ib)|{\cal P}_2=\ad\,ih|{\cal P}_2$. Since
$\|b\|={\rm Tr}|T|<\epsilon/2$, we have that
  $$
  \|\ad\,ih(x_{t_0})\|<2\epsilon.
  $$
If $\epsilon$ is sufficiently small, we have by 4.1 of
\cite{Kis99} a continuous path $(y_t)$ of unitaries in ${\cal
P}_2$ such that
  $$
  \ad\,ih(y_t)\approx0.
  $$
Since $\|\delta_{\alpha}(y_t)\|\leq \|\ad\,ih(y_t)\|+\epsilon$,
this completes the proof.
\end{pf}

\begin{lem}\label{D1}
Let $S$ be a compact infinite subset of $\R$ and $\nu$ a
probability measure on $S$ with support $S$. Let $H$ be the
multiplication operator by the identity function $x\mapsto x$ on
$L^2(\nu)$. If $\xi\in L^2(\nu)$ has norm one, there exists a
continuous path $(\xi_t)_{0\leq t<1}$ of unit vectors in
$L^2(\nu)$ such that $\xi_0=\xi$,  $(H\xi_t|\xi_t)$ and
$\|H\xi_t-(H\xi_t|\xi_t)\xi_t\|=
(\|H\xi_t\|^2-|(H\xi_t|\xi_t)|^2)^{1/2}$
are constant in $t$, and $\supp\,\xi_t$ shrinks to a three-point
set as $t\ra1$. \end{lem}
\begin{pf}
Since  both $(H\xi|\xi)=\int_Sx|\xi(x)|^2d\nu$ and
$$
\|H\xi-(H\xi|\xi)\xi\|^2=\int_Sx^2|\xi(x)|^2d\nu
-(\int_Sx|\xi(x)|d\nu)^2
$$
depend only on the modulus $|\xi(x)|$, we first choose a
continuous path $(\xi_t)_{0\leq t\leq 1}$ of unit vectors in
$L^2(\nu)$ such that $\xi_0=\xi,\ |\xi_t(x)|=|\xi(x)|$, and
$\xi_1(x)=|\xi(x)|$. Thus we may suppose that $\xi(x)\geq0$.

Let $a=\min S$, $b=\max S$, and
  \BE
  c&=&\int_S x\xi(x)^2d\nu(x),\\
  v&=&\int_S x^2\xi(x)^2d\nu(x)-c^2,
  \EE
where $c$ is the {\em mean} of $x$ and $v$ is the {\em variance}
of $x$ with respect to the probability measure $\xi(x)^2d\nu$.
Then it follows that $a<c<b$ and $0<v<(b-c)(c-a)$. (Note that a
probability measure $d\mu$ on $[a,b]$ with $\int xd\mu=c$ can be
approximated by a discrete measure
  $$
  \sum_i\lambda_i(\frac{t_i-c}{t_i+s_i}\delta_{s_i}+
  \frac{c-s_i}{t_i+s_i}\delta_{t_i}),
  $$
where $\lambda_i>0,\ \sum_i\lambda_i=1$ and $a<s_i<c<t_i<b$, whose
mean is $c$ and whose variance is
$\sum_i\lambda_i(t_i-c)(c-s_i)$). We find three distinct points
$s_1,s_2,s_3$ in $S$ such that the convex set
$$
\{\sum_{i=1}^3\lambda_i\delta_{s_i}\ |\ \lambda_i>0,\
\sum\lambda_i=1\}
$$
contains a probability measure with mean $c$ and variance $v$.
(For example, if $c\in S$ we may take $s_1=a,s_2=c,s_3=b$;
otherwise set $t_1=\max \{s\in S\ |\ s<c\}$ and
  $t_2=\min\{s\in S\ |\ s>c\}$.
Then there are three of the four points $a,t_1,t_2,b$ satisfying
the requirement. If $(b-c)(c-t_1)<v$, we may set $s_1=a,s_2=t_1,
s_3=b$; otherwise if $(t_2-c)(c-a)<v$ we may set
$s_1=a,s_2=t_2,s_3=b$; otherwise we may set
$s_1=a,s_2=t_1,s_3=t_2$.) Then for any $\epsilon>0$ we can find a
positive measurable function $g$ on $S$ with $\supp\,g\subset
\cup_i(s_i-\epsilon,s_i+\epsilon)\cap S$ such that
  \BE
  && \int g(x)d\nu=1,\\
  &&\int xg(x)d\nu=c,\\
  &&\int x^2g(x)d\nu=v+c^2.
  \EE
Define $\xi_t\in L^2(\nu)$ by
  $$
   \xi_t(x)=((1-t)\xi(x)^2+tg(x))^{1/2}.
  $$
Then $(\xi_t)_{0\leq t\leq 1}$ defines a continuous path of unit
vectors in $L^2(\nu)$ from $\xi$ to $\sqrt{g}$ such that
  \BE
  &&(H\xi_t|\xi_t)=(H\xi|\xi),\\
  &&\|H\xi_t-(H\xi_t|\xi_t)\xi_t\|=\|H\xi-(H\xi|\xi)\xi\|,
  \EE
  and $\supp(\xi_1)\subset \cup_i (s_i-\epsilon,s_i+\epsilon)\cap
S$. Continuing this argument with $\xi=\sqrt{g}$ and a smaller
$\epsilon$, we will eventually obtain a continuous path
$(\xi_t)_{0\leq t<1}$ with the required properties such that
$\cap_t\ol{\cup_{s>t}\supp\,\xi_s}=\{s_1,s_2,s_3\}$. This
completes the proof.
\end{pf}

\medskip
\small
\end{document}